\documentclass[12pt,twoside]{article}

\usepackage{amssymb,amsthm,amsmath, amsfonts}
\usepackage{graphicx}
\usepackage{color}

\renewcommand{\a}{\alpha}
\renewcommand{\b}{\beta}
\begin{document}

\centerline{}

\centerline{}

\centerline {\Large{\bf Qualitative Properties of Solutions of }}
\par \bigskip
\centerline{\Large{\bf Semilinear Elliptic Systems}}

\centerline{} \centerline{\bf {Marta
Calanchi\footnote{Universit\`a degli Studi di Milano, Italy; 
e-mail: marta.calanchi@unimi.it}, Bernhard Ruf
\footnote{Istituto Lombardo Accademia di Scienze e Lettere; e-mail: bruf001@gmail.com }}}

\centerline{}

\newtheorem{Theorem}{\quad Theorem}[section]

\newtheorem{Definition}[Theorem]{\quad Definition}

\newtheorem{Corollary}[Theorem]{\quad Corollary}

\newtheorem{Lemma}[Theorem]{\quad Lemma}

\newtheorem{Example}[Theorem]{\quad Example}
\newtheorem{Proposition}[Theorem]{\quad Proposition}

\newtheorem{Remark}[Theorem]{\quad Remark}

\begin{abstract}
The article explores the qualitative properties of solutions to elliptic equations and systems, focusing particularly on whether solutions retain the symmetry of their domains. According to the well-known Gidas-Ni-Nirenberg theorem, positive solutions to certain autonomous elliptic equations in radial domains are radial themselves. However, this symmetry can be broken in equations with power weight terms. The article also examines related results for systems of these weighted equations.
\end{abstract}
\vspace{0.76cm}
{\bf Mathematics Subject Classification :} 35J47 · 35J60 · 35J20\\
{\bf Keywords:} symmetry breaking, minimal energy solution, elliptic systems

\section{Introduction}
This article is concerned with semilinear elliptic equations and systems of such equations on symmetric domains. It is a natural question to ask whether solutions of these equations inherit the symmetry of the domain.
\par \smallskip
Indeed, a famous result of Gidas-Ni-Nirenberg (\cite{GNN}, 1979)
says that the positive solutions of certain semilinear elliptic
equations on radial domains are themselves radial. The proof of
this result relies on the so-called {\it method of moving planes}.
Generalizations of this result to more general symmetries
(symmetry with respect to hyperplanes, symmetry with respect to a
finite group of rotations) have been obtained by Pacella
\cite{pacella},  \cite{pacella1}.
\newpage
On the other hand, Smets-Su-Willem (\cite{ssw}, 2002) showed that
in related non-homogenous elliptic equations, namely equations in
which the nonlinearity has a  H\'enon type weight, the symmetry of
the positive solutions may be broken. In this article we will
concentrate on this phenomenon, in various contexts.

In the {\it first part} of this article we will describe the mentioned result of Smets-Su-Willem and related results. The methods in these results are variational, and rely on the observation that the positive solutions are often {\it minimal energy
solutions}. Indeed, Smets-Su-Willem show that for certain H\'enon weights the minimal energy solution cannot be radial, since the global minimal energy is at a lower level than the minimal energy on the subspace of radial functions.

Another possible approach to symmetry breaking is the theory of bifurcation, but we will not pursue this approach here.

 In the {\it second part} of this article, we consider these questions for coupled elliptic systems. Starting with a result by
Mitidieri (\cite{mitidieri}, 1993) we discuss the importance of
the so-called critical hyperbola for the solvability of such
systems. Then, we review a result of X.J. Wang (\cite{wang}, 1993)
concerning the application of the moving plane technique for
systems.

To discuss the critical range for systems with weights, we 
then present the Caffarelli-Kohn-Nirenberg (\cite{CKN}, 1984)
inequality and its generalization to higher order operators by
C.S. Lin (\cite{Lin}, 1986). This yields a modified critical
hyperbola which depends explicitely on the exponents of the
weights.

Finally, we turn to symmetry breaking results. We discuss
results of Cao-Peng (\cite{caopeng}, 2003) for elliptic equations
with H\'enon type weights and with nonlinearities approaching
criticality. Furthermore, we present a result
(\cite{CR}, 2010) of the authors concerning the existence of
radial and non-radial minimal energy solutions for elliptic
systems with Hardy and H\'enon type weights. In particular, we
discuss the importance of certain new hyperbolas in
describing the symmetry breaking (this result has been later  improved by Bonheure-Dos Santos-Ramos in \cite{BDR}).

\section{Symmetry and Symmetry-Breaking}
\subsection{Radial solutions}
Consider  the nonlinear elliptic equation
\begin{equation} \label{pb0}
  \begin{cases}
    -\Delta u=f(u) &\text{in $\Omega$},\\
 u>0
    &\text{in $\Omega$},\\ u=0 &\text{on $\partial \Omega$},
  \end{cases}
\end{equation}
where $f:\mathbb R\to\mathbb R$ is a continuous function,
$\Omega $ is a bounded smooth domain in $\mathbb R^N$ $(N\ge 3)$ with smooth boundary
$\partial\Omega$. We consider the functional
$$
J(u)=\frac{1}{2}\int_{\Omega}|\nabla u|^2\ dx-\int_{\Omega}F(u) \
dx
$$
associated to (\ref{pb0}). Here $F(t)$ is the primitive of $f$
which vanishes in the origin.

\
It is well known that one can obtain solutions of  (\ref{pb0}) by finding  critical points of
$J$ in the Sobolev space $H^1_0(\Omega)$,  i.e. the space of
functions whose weak derivative belongs to the space
$L^2(\Omega$). In order to have the term $\int_\Omega F(u)dx $ well-defined in $H^1_0(\Omega)$,
 one needs to impose a growth condition on $f(s)$, that is
$$
|f(s)|\le c(1+|s|^{2^*-1}), \quad s\in \mathbb R, \
2^*=\frac{2N}{N-2}
$$
so that $F$ satisfies
$$
|F(s)|\le c(1+|s|^{2^*}).
$$
Then, by the Sobolev embedding theorem $H_0^1(\Omega)\hookrightarrow
L^{2^*}(\Omega)$, $J$ is well defined on $H_0^1(\Omega)$. If we
impose the stronger growth condition
$$
|f(s)|\le c_1+c_2|s|^p, \quad 2<p+1<2^*,
$$
then by Rellich's theorem  one has a compact embedding and
consequently  solvability of the problem.
\par \smallskip
Now we concentrate on the case that $\Omega$ is a ball.  In the
celebrated paper of Gidas, Ni and Nirenberg, the authors prove the radial
symmetry of positive solutions of elliptic equation

\begin{Theorem}\label{GNN}(B. Gidas, W.M. Ni and L. Nirenberg, (1979) \cite{GNN})  Let $\Omega=\{|x|<R\}$ a ball
in $\mathbb R^N$ and $f$ be of class $C^1$. If $u>0$ is a
classical positive solution of (\ref{pb0}), then $u$ is  radially
symmetric and
\begin{equation*}
\frac{\partial u}{\partial  r}< 0, \;{\rm for}\;0<r<R.
\end{equation*}
\end{Theorem}
\par \bigskip
They also prove a similar result for  functions which depend on $r$, i.e. $f=f(r,u)$, provided that they are
decreasing in $r$:

\begin{Theorem}\label{GNN2}(B. Gidas, W.M. Ni and L. Nirenberg, (1979) \cite{GNN})
Theorem \ref{GNN} holds for $f=f(r,u)$ depending also on $r$, with
$f$ and $f_u$ continuous, provided that $f$ is {\bf decreasing} in $r$.
\end{Theorem}
As an example, consider the following equation with a {\it Hardy type}  weight:
$$
- \Delta u = \frac 1 {|x|^\beta}\, u^p \ , \ \hbox{ with } \ 0 < \beta < 2 \ .
$$
Since $\frac 1 {|x|^\beta}$  is continuous and decreasing for
$|x|\neq 0$,  the result of Gidas-Ni-Nirenberg says that the
positive solutions of this equation in the ball, and with
Dirichlet boundary conditions, are radial.

 \vspace{1cm}

 It is now natural to ask what happens if one adds an increasing radial weight?
 Let us consider for instance the {\it H\'enon equation}
\begin{equation}\label{eq:henon}
-\Delta u=|x|^{\alpha}\,u^p \ ,
\end{equation}
which was introduced in 1973 by M. H\'enon  in the context of the simulation of spherically
symmetric clusters of stars.

This equation has become a good model for a  series of problems of
great interest from a mathematical point of view, especially in the
domain of variational methods and nonlinear analysis.
In particular, much attention has been devoted to the Dirichlet problem for
positive solutions in the unitary ball $\Omega=B_1(0)$
\begin{equation} \label{pb1}
  \begin{cases}
    -\Delta u=|x|^{\alpha}\,u^p &\text{in $\Omega$},\\
 u>0
    &\text{in $\Omega$},\\ u=0 &\text{on $\partial \Omega$},
  \end{cases}
\end{equation}
for which existence, non existence, multiplicity and qualitative
properties such as radial symmetry have been considered in various
papers, under different conditions on $\alpha$ and $p$.

\par \medskip
Concerning existence results for problem
(\ref{pb1}), one first observes
that the exponent $\alpha$ affects the range of powers $p$ for which (\ref{pb1})
may possess solutions. Indeed the presence of the term $|x|^{\alpha}$ modifies the global
homogeneity of the equation and shifts upward the threshold between
existence and nonexistence (given by an application of the
Pohozaev identity).
\par \medskip
The existence of solutions in the Sobolev space $H^1_0(\Omega)$,
for $N\ge 3$, was first considered by Ni (1982) who realized that
the presence of the weight extends the range of $p$ for which a
solutions exists.

\begin{Theorem}(W.M. Ni, (1982) \cite{ni}). Problem (\ref{pb1}) possesses a {\rm radial} solution for all $p+1\in
(2,2^*_{\alpha})$, where $2^*_{\alpha}=2^*+\frac{2\alpha}{N-2}$.
\end{Theorem}

{\it Idea of proof.} The first important
ingredient is an estimate for the radial functions:
\par \smallskip \noindent
{\it {\bf Radial Lemma}. Let $u$ be a radially symmetric function on
$\Omega= B_1(0)$ with $u(1)=0$. Then,
\begin{equation}\label{radial}
|u(x)|\le\frac{||\nabla u||_{L^2(\Omega)}}{\sqrt{\omega_N(N-2)} \
|x|^{(N-2)/2}}
\end{equation}}
where $\omega_N$ is the surface area of the unit ball in $\mathbb
R^N$.
\par \smallskip \noindent
This inequality allows one to prove that the map $u\to
|x|^{\alpha/(p+1)}\cdot u$ from $H^1_{0, \mathrm{rad}}(\Omega)$ to
$L^{p+1}(\Omega)$ is compact for all $p+1 < 2^*_{\alpha}$  (here
$H_{0, \mathrm{rad}}^1(\Omega)$ denotes the space of radial
functions on $H_{0}^1(\Omega)$). Then, in order to conclude,  it
suffices to make  use of the Mountain-Pass Lemma due to Ambrosetti
and Rabinowitz \cite{AR} to obtain a solution.

\rightline{$\square$}

\par \smallskip 
Ni's result concerns a wider class of nonlinearities $f(|x|,u)$
with subcritical growth. Moreover, via a Pohozaev identity one can
prove that there exist {\it no} (nontrivial) solutions if $p+1 \ge
2^*+\frac{2\alpha}{N-2}$. We emphasize that Ni's result concerns
the existence of radial solutions. If we look for solutions in the
whole space $H_0^1(\Omega)$ we are forced to restrict ourselves to
the subcritical case $p + 1 <2^*$.
\par \medskip
We now return to the question mentioned above.  Since the function
$|x|^{\alpha}$, $\alpha > 0$, is increasing, the classical results
(moving plane techniques, maximum principle) of Gidas-Ni
-Nirenberg cannot be applied to force radial symmetry of the
solutions.

And indeed, non radial solutions appear in a natural way.

\newpage
\subsection{Symmetry breaking and multiplicity}
In an interesting paper, Smets, Su and Willem \cite{ssw} proved
some symmetry--breaking results for the H\'enon equation. Their
study was principally devoted to the symmetry of the
 so-called \emph{ground-state
  solutions}, or least energy solutions of (\ref{pb1}).

  Consider the Rayleigh quotient
\begin{equation}\label{ray}
R(u) := \frac{\int_{\Omega} |\nabla u|^2\ dx}{\left( \int_{\Omega}
|x|^\alpha |u|^{p+1} \ dx\right)^{2/(p+1)}}.
\end{equation}
For $\alpha\ge 0$ and $p+ 1 < 2^*$, it is easy to verify that
\begin{equation}\label{Sap}
S_{\alpha p}:=\inf_{\substack{u \in H_0^1(\Omega) \\ u \neq 0}}
R(u)
\end{equation}
is achieved by a function  $u$ (the ground state) which -- up to
rescaling -- is a solution of \eqref{eq:henon}.

 By the symmetry of the problem it is natural to consider also
\begin{equation} \label{eq:Saprad}
S_{\alpha p}^{rad}:=\inf_{\substack{u \in
H_{0,\mathrm{rad}}^1(\Omega)
\\ u \neq 0}} R(u)
\end{equation}
where $H_{0, \mathrm{rad}}^1(\Omega)$ denotes the space of radial
functions in $H_{0}^1(\Omega)$.

It is also easy to verify that $S_{\alpha p}^{rad}$ is achieved by
a positive function $v$, which -- after rescalings -- is a solution of
\eqref{pb1} (see Palais' symmetric criticality principle \cite{PT}). Since
the weight $|x|^{\alpha}$  is increasing, symmetrization  in
general does not increase the value of the denominator.  To the contrary,
numerical computations suggested that in some cases
\begin{equation}\label{S}
S_{\alpha p}^{rad}>S_{\alpha p} \  .
\end{equation}
Indeed, in \cite{ssw}  the authors proved that \eqref{S} holds for $\alpha$ sufficiently large,
and thus the ground state solution  cannot be radial, at least for
$\alpha$ large enough. Note that, as a consequence, \eqref{eq:henon} has at
least two solutions when $\alpha$ is large, namely a radial and a non-radial solution.

\begin{Theorem}(D. Smets, J. Su and M. Willem (2001) \cite{ssw})
\par \noindent
 Assume $N \ge 2$.  For any $2 < p+1 < 2^*$ there exists $\alpha^*>0$
such that no minimizer of R is radial provided $\alpha>\alpha^*.$
\end{Theorem}

The authors give two different proofs. We give the idea of one of
them.

\begin{proof} By a change of variable one can prove the following
asymptotic estimates for $S_{\alpha p}^{rad}$: there exists $C$
depending on $N$ and $p$ such that
$$
S_{\alpha p}^{rad}\sim C \
\left(\frac{\alpha+N}{N}\right)^{1+2/(p+1)},\ {\rm as}\ \alpha\to
+\infty.
$$
On the other hand, by testing the Rayleigh quotient $R$ on functions $u_{\alpha}$  which
concentrate near the boundary, one obtains the estimate
$$
S_{\alpha p}\le R(u_{\alpha})\le C\ \alpha^{2-N+\frac{2N}{p+1}}.
$$
Therefore
$$
S_{\alpha p}^{rad} > S_{\alpha p} \ , \ \hbox{ for } \ \alpha \ \hbox{ large } \ ,
$$
since $1+\frac 2 {p+1} > {2-N+\frac{2N}{p+1}}$ whenever $p>1$.

\end{proof}

Later on, Serra proved in \cite{serra03} the existence of at least
one non--radial solution to \eqref{eq:henon} in the critical case
$p+1=2^*$, and in \cite{badialeserra} Badiale and Serra proved the
existence of more than one solution to the same equation also for
some supercritical values of $p$. These solutions are non-radial
but have some partial symmetry, and are obtained by minimization under suitable symmetry
constraints.
\begin{Theorem}\label{badiale} (M. Badiale, E. Serra, (2004),
\cite{badialeserra}) Assume $N\ge 4$ and
$p\in(1,\frac{N+1}{N-3})$. Then, for $\alpha$ large, problem
(\ref{pb1}) has at least $\left[\frac{N}{2}\right]-1$ non radial
solutions.
\end{Theorem}

More precisely they look for solutions by minimizing the Rayleigh
quotient $R(u)$ on spaces with particular kind of symmetry: for
$N\ge 4$ and $0\le N-l\ge l$ they write $x=(y,z)\in\mathbb
R^l\times\mathbb R^{N-l}$ and consider functions $u$ in
$H^1_0(\Omega)$ which are radial with respect to $y$ and $z$ (with
an abuse of notation $u(x)=u(|y|,|z|)$). The quotient $R(u)$
attains is infimum on these spaces. Then, by considering different
$l$'s, they can distinguish these critical points one from the
other and each from the radial one.

\par \medskip
In a different approach, Cao and Peng proved in \cite{caopeng} that
for $\alpha > 0 $ fixed and for $p+1$ sufficiently close to $2^*$, the ground-state solutions of
\eqref{eq:henon} possess a unique maximum point whose distance
from $\partial \Omega$ tends to zero as $p+1 \to 2^*$; this implies again that the
ground state solutions cannot be radially symmetric.

\begin{Theorem} (D. Cao and S. Peng, (2002) \cite{caopeng})
Suppose $2<p+1<2^*$ and $\alpha>0$. Then the ground state solution
$u_p$    (after passing to a subsequence)
 satisfies that for some
$x_0\in\partial\Omega$,
\begin{enumerate}
\item[i)] $|\nabla u_p|^2\to\mu \delta_{x_0}$ as $p+1\to 2^{*}$ in
the sense of measure;

\item[ii)] $| u_p|^{2^*}\to\nu \delta_{x_0}$  as $p+1\to 2^{*}$ in
the sense of measure,

where   $\mu>0,\;\nu>0$ satisfy $\mu\ge S\nu^{\frac{2}{2^*}}$, $S$
is the best Sobolev constant, and $\delta_x$ is the Dirac mass at the point
$x$.
\end{enumerate}
\end{Theorem}

To prove their result, Cao and Peng show that the ground state solutions $u_p$
form a minimizing sequence of the best Sobolev constant $S$ as $p + 1 \to 2^*$,
and then prove the Theorem by using the concentration compactness principle
by P.L. Lions \cite{PLL}.
\par \medskip
This kind of result was improved by Peng  in \cite{peng}, where multibump
solutions for the H\'{e}non equation with almost critical Sobolev
exponent $p$ are found. The proof relies on  a finite--dimensional
Lyapunov-Schmidt reduction.  These solutions are not radial, but they are
invariant under the action of suitable subgroups of $O(N)$, and
they concentrate at boundary points of the unit ball of $\mathbb
R^N$ as $p+1 \to 2^*$. Also in this result the value of $\alpha$ remans fixed.
\par \medskip
 We also mention that  for the subcritical case $p +1 <2^*$  Byeong and
Wang studied in \cite{byeonwang} the profile of the ground state
solutions, proving thatthey concentrate
on the boundary, for $\alpha\to +\infty$.
\par \medskip
Later on, Cao, Peng and Yan \cite{caopengyan} improved this
result. They proved that for fixed $p +1 <2^*$ the maximum point
$x_{\alpha}$ of the ground state solution $u_{\alpha}$ satisfies
$\alpha(1 - |x_{\alpha}|) \to l\in (0,+\infty)$ as
$\alpha\to+\infty$. They also obtain
 the existence of multi-peaked solutions and give their asymptotic behavior.
\par \smallskip
  We also mention a very recent result, due to Y. Zhang and J. Hao \cite{ZH}, about
the H\'enon biharmonic equation
\begin{equation}\label{biharmonic}
\begin{cases}
   \Delta^2 u= |x|^{\alpha}\,u^p &\text{in $\Omega$},\\
 u>0
 &\text{in $\Omega$},\\ u=\frac{\partial u}{\partial \nu}=0 &\text{on $\partial \Omega$},
  \end{cases}
\end{equation}
where $\nu$ is the unit exterior normal on $\partial\Omega$,
$\alpha>0$ and $2< p+1 < 2^{**}=\frac{2N}{N-4}$, where $2^{**}$ is
the critical Sobolev exponent for the embedding
$H^2_0(\Omega)\hookrightarrow L^{p+1}(\Omega)$. The authors are
interested in the asymptotic behavior of the ground-state
solutions of (\ref{biharmonic}). More precisely, define
$$
S_{\alpha, p}=\inf_{\substack{u \in H_0^2(\Omega) \\ u \neq
0}}\frac{\int_{\Omega}|\Delta
u|^2}{\left(\int_{\Omega}|x|^{\alpha}|u|^{p+1}\right)^{2/(p+1)}}.
$$
Since the embedding $H^2_0(\Omega)\hookrightarrow L^{p+1}(\Omega)$ is
compact for $p+1<2^{**}$, there exists a nontrivial function $u_p\in
H^2_0(\Omega)$ such that $S_{\alpha, p}$ is achieved by $u_p$.
In analogy to the result of Cao and Peng \cite{caopeng} they prove

\begin{Theorem}(Y. Zhang and J. Hao, 2011, \cite{ZH}).
Suppose $2<p+1 <2^{**}$ and $\alpha>0$. Then the ground state
solutions $u_p$    (after passing to a subsequence)
 satisfy that for some
$x_0\in\partial\Omega$,
\begin{enumerate}
\item[i)] $|\Delta u_p|^2\to\mu \delta_{x_0}$ as $p+1\to 2^{**}$ in
the sense of measure;

\item[ii)] $| u_p|^{2^{**}}\to\nu \delta_{x_0}$ as $p+1\to 2^{**}$
in the sense of measure,

where   $\mu>0,\;\nu>0$ satisfy $\mu\ge {\bar
S}\nu^{\frac{2}{2^*}}$, $\delta_x$ is the Dirac mass at $x$, and
$\bar S$ is the best Sobolev constant
$$
{\bar S}=\inf_{\substack{u \in H_0^2(\Omega) \\ u \neq
0}}\frac{\int_{\Omega}|\Delta
u|^2}{\left(\int_{\Omega}|u|^{2^{**}}\right)^{2/{2^{**}}}}.
$$
\item[iii)] when $p+1$ is close to $2^{**}$,  $u_p$ has a unique
maximum point $x_p$ such that $dist(x_p,\partial\Omega)\to 0$ as
$p+1\to 2^{**}$.
\end{enumerate}
\end{Theorem}
Thus, also for the biharmonic equation we find symmetry breaking: the
ground state solutions are non radial if $p+1$ is sufficiently close to $2^{**}$.

\section{Symmetry Breaking in Systems}
\subsection{The critical hyperbolas}
Let us now consider the following system of coupled elliptic equations
with different H\'enon coefficients:
\begin{equation} \label{pb2}
  \begin{cases}
    -\Delta v=|x|^{\alpha}u^{p} &\text{in $\Omega$},\\
-\Delta u=|x|^{\beta}v^{q} &\text{in $\Omega$},\\ u>0 \ , \
 v>0 &\text{in $\Omega$},\\ u=v=0 &\text{on $\partial \Omega$},
  \end{cases}
\end{equation}
where $\Omega$ is a ball in $\mathbb R^N$ centered in the origin,
$N\ge 3$, $p,q>1$, and $\alpha, \beta>-N$. In particular, we
discuss existence, multiplicity and qualitative properties (such
as radial symmetry in the case $\Omega$ a ball) of solutions.

 The case
$\alpha=\beta=0$  has been studied by many authors.
 Here the natural restriction on the exponents $p$ and $q$ for existence\,/\,non-existence of solutions is given
by the so-called {\it critical hyperbola}, that is
\begin{equation}\label{crit-hyp0}
\frac{N}{p+1}+\frac{N}{q+1} = N-2 \ .
\end{equation}
This hyperbola was first introduced by E. Mitidieri
\cite{mitidieri} who proved non-existence of solutions for
$(p+1,q+1)$ lying on or above the hyperbola, using a Pohozaev type
identity. Existence of solutions for $(p+1,q+1)$ below the
critical hyperbola was proved by de Figueiredo-Felmer (see
\cite{deF-F}) and by Hulshoff-van der Vorst (see \cite{H-vdV},
\cite{hulsvander2}) by using a variational set-up with fractional
Sobolev spaces. A different approach, working with Sobolev-Orlicz
spaces (which allows a generalization to non-polynomial
nonlinearities) can be found in  \cite{Clem}, \cite{deF-doO-R}.

The general case $\alpha \neq 0$, and/or $\beta\neq
0$ was investigated independently by de\,Figueiredo-Peral-Rossi \cite{DPR} and Liu Fang-Yang Jianfu
\cite{liuyang}; in both papers an approach via fractional Sobolev
spaces is used. As in the scalar case, the presence of the weight
functions $|x|^{\alpha}$ and $|x|^{\beta}$ affects the range of
$p$ and $q$ for which the problem may have solutions. Indeed, in
\cite{DPR} and \cite{liuyang} it is shown that the dividing line
between existence and non-existence is given by  the following
``weighted" critical hyperbola
\begin{equation}\label{crit-hyp2}
\frac{N +\alpha}{p+1}+\frac{N+\beta}{q+1} = N-2.
\end{equation}

For future reference, we call the hyperbola (\ref{crit-hyp0}) the
M-hyperbola (for {\it Mitidieri hyperbola}), and the hyperbola
(\ref{crit-hyp2}) the $\alpha\beta$-hyperbola).

\begin{Theorem} (D. de Figueiredo, I. Peral, J.D.  Rossi, \cite{DPR})
\par \noindent
Assume that $p,\ q > 1$ and $-N<\alpha,\ \beta < 0$ verify
\begin{equation}\label{crit-hyp4}
\frac{N +\alpha}{p+1}+\frac{N+\beta}{q+1} >N-2.
\end{equation}
\begin{equation*}
\frac{1}{p+1}+\frac{1}{q+1} < 1.
\end{equation*}
and
\begin{equation}
p+1<\frac{2(N+\alpha)}{N-4} \;\; and \;\;
q+1<\frac{2(N+\beta)}{N-4}\;\; if \,\, N\ge 5.
\end{equation}
Then there exists at least one positive solution of  (\ref{pb2}).
\end{Theorem}

For plots of the $\alpha\beta$-hyperbola, see figures 1 - 3 below.
\par \bigskip

\begin{Theorem} (F. Liu and J. Yang, (2007) \cite{liuyang})
\par \noindent
Let  \ $\alpha,\ \beta>-N$ and $ p,\ q>1$ such that
\begin{equation*}
\frac{N +\alpha}{p+1}+\frac{N+\beta}{q+1} >N-2.
\end{equation*}
If one of the following assumptions is satisfied
\begin{enumerate}
\item[i) ] (Hardy-case) \ $0\ge\alpha,\ \beta>-N$,
$\alpha+\beta>-4$;

\item[ii) ] (H\'enon-case)\ $\alpha, \beta\ge 0$, $\frac{\alpha
N}{N+\alpha}+\frac{\beta N}{N+\beta}<4$;

\item[iii) ] (Mixed-case)\ $-N<\beta\le 0\le \alpha$,
\end{enumerate}
then there exists at
least one positive solution of (\ref{pb2}).
\end{Theorem}

The authors proved this result for more general nonlinearities with
some growth conditions.
\par \medskip \noindent
{\it Idea of Proof:} In both papers  the following functional associated to the system
is considered
\begin{equation}\label{J}
J(u,v) = \int_\Omega \nabla u \nabla v - \frac 1{p+1}\int_\Omega |x|^\alpha |u|^{p+1} -
\frac 1{q+1}\int_\Omega |x|^\beta |u|^{q+1}
\end{equation}
 A natural choice of the
space to define the functional would seem to be $H_0^1(\Omega)
\times H_0^1(\Omega)$; however, this is too restrictive and does
not give an optimal result. The idea is to ``destroy'' the
symmetry between $u$ and $v$ by requiring more regularity of $u$
than of $v$ if $p$ is large and $q$ is small (and viceversa). 
More
precisely, the authors consider {\it fractional Sobolev spaces}
$H^s(\Omega)$, which can be defined via interpolation, or
 by Fourier series: let $(e_i)_{i \in \mathbb N}$ be the orthogonal basis of eigenfunctions of the Laplacian in $H_0^1(\Omega)$, and $\lambda_i$ the corresponding eigenvalues. Then,  define the Hilbert space $H^s(\Omega) := \big\{ u : \sum_{i \in \mathbb N} \lambda_i^s (u,e_i)_{L^2}^2 < \infty \big\}$ of functions with fractional derivative of order $s$, and analogously, the fractional derivative operator $A^s: u \mapsto \sum_{i \in \mathbb N}  \lambda_i^{s/2}  (u,e_i)_{L^2}\,e_i$. With these definitions, one sees that the
 integral $\int_\Omega \nabla u \nabla v\,dx $ can be defined on $H^s \times H^t$, with
 $s + t = 2$, by performing a ``fractional partial integration'', that is by writing
 $$
 \int_\Omega \nabla u \nabla v \, dx = \int_\Omega A^s u A^t v \,dx : H^s \times H^t \to \mathbb R
 $$
 To obtain that $J$ is a well-defined functional on $H^s \times H^t$ one now uses the
 Sobolev embeddings of $H^s(\Omega)$ into weighted $L^p$-spaces:
 $$
 H^s(\Omega) \subset L^p(\Omega, |x|^\alpha dx) \ , \ \hbox{ for } \  p \le \frac {2(N+\alpha)}{N-2s} \ ,
 $$
with compact inclusions for strict inequality. The conditions on $s, p , \alpha$ and $t,q,\beta$, and the condition $s + t = 2$, yield the hyperbola.

A second difficulty (which is common to such systems) is the strong indefiniteness of the functional $J$. A way to get around this difficulty is to consider a {\it finite dimensional approximation}, then do a standard linking approach in the finite-dimensional spaces to obtain approximate solutions, and then prove the convergence of the sequence of approximate solutions to a solution of the system.

\rightline{$\square$}

\vspace{1cm}
\subsection{Non existence: a generalized identity \\ of Pohozaev type  }
The previous results show that below the critical $\alpha\beta$-hyperbola (\ref{crit-hyp2}) one has the existence of solutions. One expects non-existence of solutions on or above the critical hyperbola.
This can indeed be shown, by extending the argument of Mitidieri  in \cite{mitidieri} to the non-autonomous case, by deriving a generalized identity of Pohozaev type.

\par \medskip

\begin{Theorem}(M. Calanchi, B. Ruf \cite{CR}).
Let  $\Omega\subset \mathbb R^N $  be a bounded, smooth,
starshaped domain with respect to $0\in\mathbb R^N$.
 If $\alpha,\ \beta>-N$ and

\begin{equation}
\frac{N+\alpha} {p+1}+\frac{N+\beta}{q+1}\le N-2
\end{equation}
then the problem
\begin{equation} \label{pb}
  \begin{cases}
    -\Delta v=|x|^{\alpha}|u|^{p-1} u&\text{in $\Omega$},\\
-\Delta u=|x|^{\beta}|v|^{q-1}v &\text{in $\Omega$},\\  u=v=0
&\text{on $\partial \Omega$},
  \end{cases}
\end{equation}
has no nontrivial strong positive solutions.
\end{Theorem}

We give an idea of the proof in the Hardy case (for $\alpha,
\beta$ negative).

\begin{proof}
For this case we follow the idea developed by B. Xuan {(see
Appendix in  \cite{Xuan})}.

Let $(u, v)$ be a positive solution of system (\ref{pb2}). Due to
the Hardy weights
 this solution may be singular in the origin, but
standard regularity results imply that for every $\delta$ small,
$u$ and $v$ belong to $C^2(\Omega\setminus B_{\delta}(0))\cap
C^0(\overline {\Omega\setminus B_{\delta}(0)})$. We multiply the
equations respectively by $(x\cdot\nabla u)$ and by $(x\cdot\nabla
v)$, add the two equations and integrate over $\Omega_{\delta}=\Omega\setminus B_{\delta}(0)$
\begin{multline}{\label{delta}}
-\int_{\Omega_{\delta}}\{\Delta u \,(x\cdot\nabla v)+\Delta v
\,(x\cdot\nabla u)\}dx =
 \\
\int_{\Omega_{\delta}} \left\{\displaystyle\frac{\partial
G}{\partial u}(x\cdot\nabla u)+\frac{\partial G}{\partial
v}(x\cdot\nabla v)\right\}dx \ ,
\end{multline}
where
\begin{equation} \label{G} G(x,u,v)=\frac{1
}{p+1}|x|^{\alpha}|u|^{p+1}+\frac{1 }{q+1}|x|^{\beta}|v|^{q+1}
\end{equation}

Applying the Divergence Theorem to  $x G$, one has for the right
side of \eqref{delta}
\begin{multline}\label{divGdelta}
\int_{\Omega_{\delta}}\Big\{\displaystyle\frac{\partial
G}{\partial u}(x\cdot\nabla u)+\frac{\partial G}{\partial
v}(x\cdot\nabla v)\Big\}dx = \\
=
-\frac{N-|\alpha|}{p+1}\int_{\Omega_{\delta}}\frac{u^{p+1}}{|x|^{|\alpha|}}dx -
\frac{N-|\beta|} {q+1}\int_{\Omega_{\delta}}\frac{v^{q+1}}{|x|^{|\beta|}}dx + \\
 + \int_{|x|=\delta}G(u,v,x)(x\cdot n) ds\ ,
\end{multline}
while for the left side of  \eqref{delta} (see \cite{mitidieri},
Corollary 2.1) one has
\begin{multline}
\int_{\Omega_{\delta}}\Big\{\Delta u \,(x\cdot\nabla v)+\Delta v
\,(x\cdot\nabla u)\Big\}dx =
\\
=\int_{\partial\Omega_{\delta}}\Big\{\frac{\partial u}{\partial
n}(x\cdot \nabla v)+\frac{\partial v}{\partial n}(x\cdot \nabla
u)-(\nabla u\nabla v)(x\cdot n)\Big\}ds
+(N-2)\int_{\Omega_{\delta}}(\nabla u\nabla v)dx
 \\
=\int_{|x|={\delta}}\Big\{\frac{\partial u}{\partial n}(x\cdot
\nabla v)+\frac{\partial v}{\partial n}(x\cdot \nabla u)-(\nabla
u\nabla v)(x\cdot n)\Big\}ds
\\
+\int_{\partial\Omega}\Big\{\frac{\partial u}{\partial
n}\frac{\partial v}{\partial n}(x\cdot n)\Big\}ds
+(N-2)\int_{\Omega_{\delta}}(\nabla u\nabla v)dx.
\end{multline}
Now, multiplying the first equation by $u$, the second by $v$ and
integrating, one obtains
$$
\int_{\Omega_{\delta}}(\nabla u\nabla v)dx -
\int_{|x|={\delta}}v\nabla u\cdot x =
\int_{\Omega_{\delta}}|x|^{\beta}|v|^{q+1}dx
$$
and
$$
\int_{\Omega_{\delta}}(\nabla u\nabla v)dx -
\int_{|x|={\delta}}u\nabla v\cdot x =
\int_{\Omega_{\delta}}|x|^{\alpha}u^{p+1}dx
$$
Then one proves that  all the integrals along $\{ |x|=\delta\}$ go
to zero, at least for a subsequence $\delta_k\to 0$. So that one
has  the { identity}

\begin{multline}\label{pohozaev}
\int_{\partial\Omega}\Big\{\frac{\partial u}{\partial
n}\frac{\partial v}{\partial n}(x\cdot n)\Big\}ds \\
=\Big\{-(N-2) + \frac{N-|\alpha|} {p+1}+\frac{N-|\beta|}
{q+1}\Big\}\int_{\Omega}|x|^{\alpha}|u|^{p+1}dx
\end{multline}
which gives non-existence of positive solutions on the critical
hyperbola.

\end{proof}

\subsection{Existence: the ground state}

In this section we are interested in a particular type of
solution: the  ground state. This needs a different approach  with
respect to the one used for instance by de Figueiredo-Felmer
(\cite{deF-F}) and by Hulshoff-van der Vorst (see \cite{H-vdV}),
 or by de Figueiredo-Peral-Rossi (\cite{DPR}) and Liu-Yang
(\cite{liuyang}) for the weighted systems.

\vspace{0.6cm}
We begin with the case without weights:
\begin{Theorem}\label{wang}(X. J. Wang, (1993) \cite{wang})
We consider the following problem
\begin{equation}\label{S_q}
S_{p,q}=\inf_{W^{2,r}_{\theta}(B)}\frac{\int_B|\Delta
u|^{r}}{\left(\int_B|u|^{p+1}\right)^{r/p+1}},\quad
r=\frac{q+1}{q}
\end{equation}
where  $B$ is the unit ball in $\mathbb R^N$, and
$W^{2,r}_{\theta}(B)$ denotes the set of functions in $W^{2,r}(B)$
which vanish on $\partial B$ and $p$ and $q$ satisfy
\begin{equation}\label{M crit-hyp}
\frac{N}{p+1}+\frac{N}{q+1} > N-2 \
\end{equation}
(i.e. $(p,q)$ lies below the M-hyperbola). Then the infimum of \eqref{S_q} is achieved by a
positive function $u$ which is spherically symmetric and radially
decreasing. Moreover, $u$ is a weak solution of the problem
\begin{equation} \label{eqwang}
  \begin{cases}
   -\Delta\left((-\Delta u)^{1/q}\right)=S_{p,q}\,u^p &\text{in $B$},\\
 u=\Delta u=0 &\text{on $\partial B$},
  \end{cases}
\end{equation}

\end{Theorem}
\par \smallskip
\begin{Remark} Set $v=(-\Delta u)^{1/q}$, then
$(u,v)$ is a weak solution of
\begin{equation} \label{eqwang2}
  \begin{cases}
    -\Delta u=v^{q} &\text{in $B$},\\
-\Delta v=S_{p,q}\, u^p &\text{in B},\\
 u=v=0 &\text{on $\partial B$},
  \end{cases}
\end{equation}

\end{Remark}

\vspace{0.6cm} Analogously  one can (formally) deduce from the
second equation in (\ref{pb2})
$$
v = (-\Delta u)^{\frac 1 q} |x|^{-\frac{\beta}q} \ ,
$$
and inserting this into the first equation one obtains the
following scalar equation for the $u$-component
\begin{equation} \label{eq2}
-\Delta\big((-\Delta
u)^{1/q}|x|^{-{\beta}/{q}}\big)=|x|^{\alpha}u^{p}
\end{equation}
It is interesting to investigate the r\^ole played by the weights
$\alpha$ and $\beta$ when dealing with the existence and symmetry
of {\it ground state} (or {\it minimal energy}) solutions of
equation (\ref{eq2}), that is, minimizers $u$ of the corresponding
Rayleigh quotient
\begin{equation}\label{Req}
\displaystyle R(u)=\frac{\int_{\Omega}|x|^{-{\beta(r-1)}}|\Delta u|^r
}{\left(\int_{\Omega}|x|^{\alpha}|u|^{p+1}\ dx
\right)^{\frac{r}{p+1}}} \quad , \quad r := \frac{q+1}q \ ,
\end{equation}
on the  weighted Sobolev space
$$
W^{2,r}(\Omega, |x|^{-{\beta}(r-1)}dx)\cap W^{1,r}_0(\Omega) \
$$
Here we denote  with $W^{2,r}(\Omega,|x|^{-{\beta}(r-1)}dx)$ the
set of functions $u\in W^{2,1}_{\rm loc}(\Omega)$ such that
$$
\int_{\Omega}(|u|^r+|\nabla u|^r+\sum_{|\xi|=2}|D^{\xi}u|^r
|x|^{-{\beta}(r-1)}) \ dx<+\infty \ ,
$$
endowed with the norm
$$
 \|u\|_{W^{2,r}(\Omega,|x|^{-{\beta}/{q}}dx )}:= \bigg(\int_{\Omega}(|u|^r+|\nabla u|^r+\sum_{|\xi|=2}|D^{\xi}u|^r
|x|^{-{\beta}(r-1)}dx) \ dx\bigg)^{1/r} \ ;
 $$
also, we denote with
 $$
W^{2,r}_{rad}(\Omega,|x|^{-{\beta}(r-1)}dx)
 $$
 the subspace of
 $W^{2,r}(\Omega,|x|^{-{\beta}(r-1)}dx)$ of radial functions.
 Furthermore, let
 $$
 {W^{2,r}_D(\Omega,|x|^{-{\beta}(r-1)}dx )}
 $$ denote the closure of
$\{\phi\in C^{\infty}(\Omega):\;\phi=0\; {\rm on}\;
\partial \Omega\}$ in ${W^{2,r}(\Omega,|x|^{-{\beta}(r-1)}dx )}$, i.e. the closure of the smooth functions in $\Omega$
with Dirichlet  boundary conditions, and with
$$
 {W^{2,r}_{D,rad}(\Omega,|x|^{-{\beta}(r-1)}dx )}
$$
the corresponding subspace of radial functions (for $\Omega$ a
ball).
\par \smallskip
Critical points of $R(u)$ on
$W^{2,r}_D(\Omega,|x|^{-{\beta}/{q}}dx)$ are (up to rescaling)
weak solutions of (\ref{eq2}), i.e. they verify
$$
\begin{cases}
\displaystyle \int_{\Omega}(-\Delta u)^{1/q} |x|^{-{\beta}/{q}}\,
(-\Delta\varphi) dx = \int_{\Omega}|x|^{\alpha}u^{p}\varphi \
dx \ , \vspace{0.2cm} \\
\displaystyle  \quad \hbox{for all }  \varphi\in
W^{2,r}_D(\Omega,|x|^{-{\beta}/{q}}dx) \ ,
\end{cases}
$$
and, moreover,  if $v = (-\Delta u)^{\frac 1 q}
|x|^{-\frac{\beta}q}$, then $v\in
W^{2,\frac{p+1}{p}}(\Omega,|x|^{-\frac{\alpha}{p}}dx)\cap
W^{1,\frac{p+1}{p}}_0(\Omega)$. In accordance, by a {\it strong
solution} of the system we mean a couple $(u,v)$ of weak solutions
such that
$$(u,v)\in
W^{2,r}(\Omega,|x|^{-\frac{\beta}{q}}dx)\cap W^{1,r}_0(\Omega)
\times W^{2,\frac{p+1}{p}}(\Omega,|x|^{-\frac{\alpha}{p}}dx)\cap
W^{1,\frac{p+1}{p}}_0(\Omega).$$

We give first the following existence result for Hardy-type
systems

\begin{Theorem}\label{existence} (Hardy case; M. Calanchi, B. Ruf (2010) \cite{CR}).
\par \smallskip \noindent
Let $ 0\ge \alpha, \beta > -N$, and  $p,q > 1$.
 If $\label{crit-hyp2b} \frac{N-|\a|}{p+1}+\frac{N-|\beta|}{q+1} > N-2$
 (i.e. $(p+1,q+1)$ lies below the $\alpha\beta$-hyperbola),  then
  system
(\ref{pb1}) has a positive nontrivial solution $(\bar u,\bar v)$.
Moreover if $\Omega$ is a ball then $(\bar u,\bar v)$ is radially
symmetric.
\end{Theorem}

\begin{figure}[h]
 \centerline{\includegraphics[width=13cm]{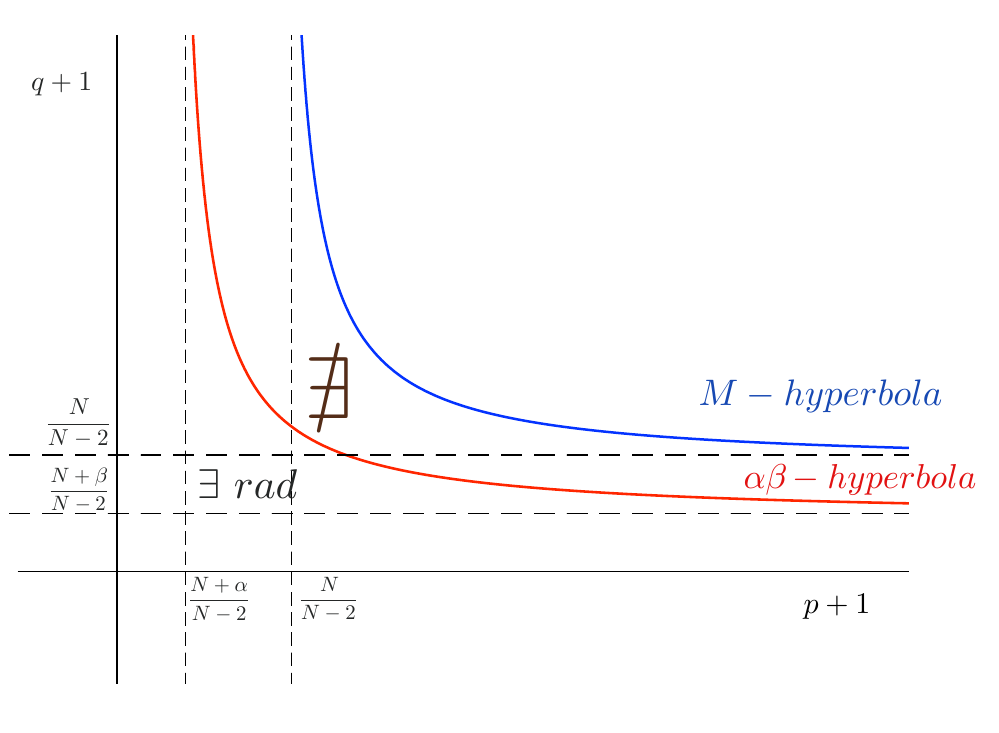}}
  \centerline{{\it Figure 1: Hardy case: $\alpha, \beta < 0$}}
  \end{figure}

\begin{proof}
 From  a generalization due to Lin (see \cite{Lin})
of the Caffarelli-Kohn-Nirenberg inequality \cite{CKN}, one can
prove that the  embedding
$$
W^{2,r}(\Omega,|x|^{-{\beta}(r-1)}dx)\hookrightarrow
L^{p+1}(\Omega,|x|^{\alpha})
$$
is compact. To prove the existence of a positive solution we
minimize the Rayleigh quotient $R(u)$ given in (\ref{Req}). By the
compactness of the embedding  the infimum is attained by a
positive function, which is a (strong) solution of problem
(\ref{pb2}).
 For $\alpha = \beta = 0$
and if $\Omega$ is a ball, as  proved by X.J. Wang \cite{wang},
 the ground state of (\ref{Req}) is a radial and radially
decreasing positive function (Theorem \ref{wang}). By adapting his
argument (moving planes technique and maxi\-mum principle) one can
extend this result also to the values $\alpha, \beta < 0$ (noting
that the weights do not interfere with the moving planes
technique, since they are decreasing).
\par \medskip
\end{proof}

\subsection{Symmetry breaking for systems}

As in the case of a single equation, there are two type of symmetry breaking
results: one for fixed $(\alpha,\beta)$ and $(p,q)$ near the critical hyperbola, and the other one for fixed $(p,q)$ with $\alpha$ large:

\begin{Theorem} (H. He, J. Yang (2008) \cite{HY})
\par \smallskip \noindent
Consider the (H\'enon) problem
\begin{equation} \label{pbepsilon}
  \begin{cases}
    -\Delta v=|x|^{\alpha}u^{p} &\text{in $\Omega$},\\
-\Delta u=|x|^{\beta}v^{q_{\varepsilon}} &\text{in $\Omega$},\\
u>0 \ , \  v>0
    &\text{in $\Omega$},\\ u=v=0 &\text{on $\partial \Omega$},
  \end{cases}
\end{equation}
 where $0\le\alpha<pN, \;\beta>0,\; N\ge 8,\; p>1,\;
 q_{\varepsilon}>1.$ Assume that $q_{\varepsilon}\to q$ as ${\varepsilon}\to
 0^+$ and $q_{\varepsilon},\; q$ satisfy respectively
\begin{equation}\label{crit-hyp-varepsilon}
\frac{N +\alpha}{p+1}+\frac{N+\beta}{q_{\varepsilon}+1} >
N-2,\quad\frac{N +\alpha}{p+1}+\frac{N+\beta}{q+1} = N-2.
\end{equation}
Then the ground state solution $(u_{\varepsilon},
v_{\varepsilon})$, up to subsequence, satisfies that for some
$x_0\in\partial\Omega$,
\begin{enumerate}
\item[i)]as ${\varepsilon}\to 0$, $|\Delta
v_{\varepsilon}|^{\frac{p+1}{p}}\rightharpoonup\mu\delta_{x_0}$,
$|x|^{\frac{\alpha(p+1)}{p}}u^{p+1}_{\varepsilon}\rightharpoonup\mu\delta_{x_0}$
in the sense of measure;
\item[ii)] $|
v_{\varepsilon}|^{{q+1}}\rightharpoonup\nu\delta_{x_0}$ in the
sense of measure,

where $\mu>0,\;\nu>0$ satisfy $\mu\ge
S_{p,q}\nu^{\frac{p+1}{p(q+1)}},\;\delta_x$ is the Dirac mass at
$x$, and
$$
S_{p,q}=\inf_{{W^{2,\frac{q+1}{q}}_{\theta}(\Omega)}\setminus
\{0\}}\frac{\int_{\Omega}|\Delta u|^{\frac{q+1}{q}}\ dx
}{\left(\int_{\Omega}|u|^{p+1}\ dx \right)^{\frac{q+1}{q(p+1)}}}
$$
is the Sobolev constant.

\end{enumerate}
\end{Theorem}
Thus, in high dimensions ($N \ge 8$), the minimal energy solutions concentrate near
a boundary point if $(p,q)$ tends to the critical hyperbola, which implies in particular that for $\Omega = B_R(0)$ the ground state solution cannot be radial.

\begin{Theorem}\label{calruf1} (H\'enon type system, M. Calanchi, B. Ruf (2010),\cite{CR}
).
\par \smallskip \noindent

Let $\alpha, \beta
> 0$, and suppose that $p,q> 1$.
\par \smallskip \noindent
a) If \ $\label{crit-hyp2c} \frac{N+\a}{p+1}+\frac{N+\beta}{q+1}
\le N-2$ (i.e., on or above the $\alpha\beta$-hyperbola) and
$\Omega$ is starshaped,
 then system (\ref{pb1}) has no  non-trivial solutions.
\par \medskip \noindent
\par \noindent
b) Suppose that $\Omega=B_1(0)$. If  $q > \max\{1,\frac \beta N\}$
and  $\frac{N+\a}{p+1}+\frac{N+\beta}{q+1} > N-2$ (i.e. below the
$\alpha\beta$-hyperbola),
 then system (\ref{pb1}) has a radial solution (not necessarily of minimal energy)
\par \medskip \noindent
c) If  $\frac{N}{p+1}+\frac{N}{q+1} > N-2$ (i.e. below the
M-hyperbola), then $\inf_{E} R(u)$ is attained and hence system
(\ref{pb1}) has a ground state solution; furthermore, if
$p>\frac{N}{(N-1)q-1}$ and $\alpha
> 0$ is sufficiently large, then the ground state solution is {\rm
not radially symmetric}.

\end{Theorem}

The situation is illustrated in the following figure 2:

\begin{figure}[h]
 \centerline{\includegraphics[width=13.5cm]{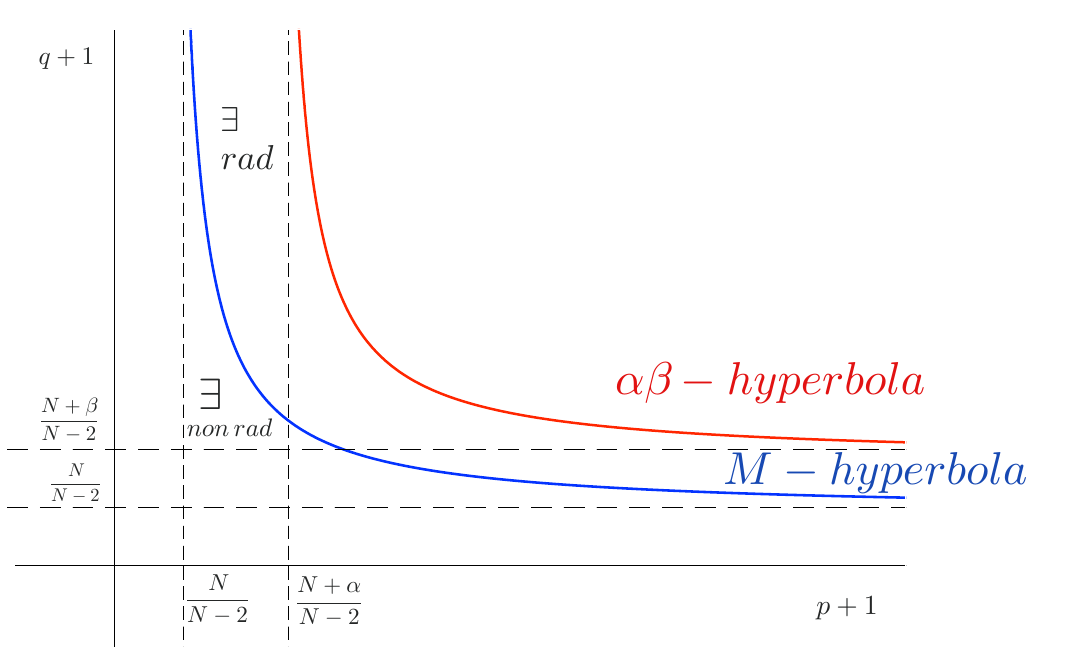}}
  \centerline{{\it Figure 2: Henon case: $\alpha, \beta > 0$}}
  \end{figure}

Since the original article contains an inaccuracy in the proof, 
we present here the essential steps.

\begin{proof} $a)$ has been already proved. 

\medskip

$b)$ For  the existence,
as in Ni's result, the  essential  ingredient is  an estimate for
the radial functions which allows one to prove that the embedding
$$
E_{rad}=W^{2,r}_{D,rad}(\Omega,|x|^{-\beta/q}dx)\hookrightarrow
L^{p+1}(\Omega,|x|^{\a}dx)\ , \ r = \frac{q+1}q
$$
is continuous and compact. Therefore, the infimum
$$
S_{\a,\beta}^{rad}=\inf_{u\in E_{rad}\setminus
\{0\}}\frac{\int_{\Omega}|x|^{-\beta/q}|\Delta u|^r\ dx
}{\left(\int_{\Omega}|x|^{\alpha}|u|^{p+1}\ dx
\right)^{\frac{r}{p+1}}}
$$
is attained by a function $\bar v$ which is, after rescaling, a
solution of (\ref{pb2}).

\medskip
$c)$ First, 
we have to show that $S_{\alpha,\beta} := \inf_E R(u)$ is
attained. Let $\{u_n\} \subset E$ be a minimizing sequence. We may
assume that
$$
\int_\Omega |x|^{\alpha}|u_n|^{p+1} = 1 \ ,\ \hbox{ and } \
\int_\Omega |x|^{-\beta/N}|\Delta u_n|^r dx \to m > 0 \ .
$$
Then clearly $\int_\Omega |\Delta u_n|^r dx \le c$, and by the
assumption and the Rellich-Kondrachov compactness theorem it
follows that $\{u_n\}$ has a convergent subsequence in
$L^{p+1}(\Omega)$, and hence also in $L^{p+1}(\Omega,|x|^\alpha
dx)$. This is sufficient to conclude that $\inf_E R(u)$ is
attained.
\par \smallskip

We show that if $\alpha > 0$ is sufficiently large, then the
radial ground state level lies above the ground state level:
indeed, by Proposition \ref{10} below, we have the following lower
estimate for the radial ground state level:
$${\label{radlev}}
S_{\a,\b}^{rad}\ \geq\ C\ \alpha^{2r+\frac{r}{p+1}-1-\frac1q}\ ,\
\hbox{ for } \ \alpha\geq\alpha_0
$$
On the other hand, for the ground state level the following upper
estimate holds:
 There exist  $C>0$ and
$\a_0$ such that for $\a\ge\a_0$
$${\label{Nradlevel}}
S_{\a,\beta}\ \le\ C\ \a^{2r-N+N\frac{r}{p+1}}
$$

From these two inequalities it follows that the ground state is
non radial for $\alpha$ sufficiently large, and for
$p>\frac{N}{(N-1)q-1}$ (for example this is true if
$p>\frac{N}{N-2}$).

%

\noindent
We give now an estimate from below for the radial level
$$
S_{\a,\beta}^{rad}=\inf_{u\in E_{rad}\setminus
\{0\}}\frac{\int_{\Omega}|x|^{-\beta/q}|\Delta u|^r\ dx
}{\left(\int_{\Omega}|x|^{\alpha}|u|^{p+1}\ dx
\right)^{\frac{r}{p+1}}}
$$
\begin{Proposition}\label{10}

  {\it  There exist $C>0$ and $\alpha_0$ such
that

$$
S_{\a,\beta}^{rad}\geq C
\alpha^{2r+\frac{r}{p+1}-1-\frac1q},\qquad\alpha\geq\alpha_0
$$
}
\end{Proposition}

\medskip
 Indeed, let $\varepsilon=\frac{N}{N+\a}$ and $u(x)=u(|x|)$ a smooth radial
 function such that $u=0$ on $\partial\Omega$. Let $v(\rho)=u(\rho^{\varepsilon})$.
 We have
 $$v'(\rho)=\varepsilon u'(\rho^{\varepsilon})\rho^{\varepsilon -1}\quad {\rm and}\quad
 v''(\rho)=\varepsilon^2 u''(\rho^{\varepsilon})\rho^{2\varepsilon-2}
 +\varepsilon(\varepsilon-1)u'(\rho^{\varepsilon})\rho^{\varepsilon-2}
 $$
so that
$$
u'(\rho^{\varepsilon})=\rho^{1-\varepsilon}\varepsilon^{-1}v'(\rho)
\quad {\rm and}\quad
u''(\rho^{\varepsilon})=\varepsilon^{-2}\rho^{2-2\varepsilon}[v''(\rho)-(\varepsilon-1)\rho^{-1}v'(\rho)]
$$
Therefore, by the change of variable $t=\rho^{\varepsilon}$,
  \[
  \int_{\Omega}|x|^{-\beta/q}|\Delta u|^r
  dx=\omega_{N-1}\int_0^1\Big|u''(t)+\frac{N-1}{t}u'(t)\Big|^rt^{N-1{-\beta/q}}dt
  \]
  \[
 =\omega_{N-1}\int_0^1\varepsilon\rho^{\varepsilon N-\frac{\varepsilon\beta}{q}-1}
 \Big|u''(\rho^{\varepsilon})+\frac{N-1}{\rho^{\varepsilon}}u'(\rho^{\varepsilon})\Big|^r d\rho
 \]
\[
 =\omega_{N-1}\int_0^1\varepsilon\rho^{\varepsilon N-\frac{\varepsilon\beta}{q}-1}\Big|\varepsilon^{-2}\rho^{2-2\varepsilon}
 \Big[v''(\rho)-(\varepsilon-1)\rho^{-1}v'(\rho)\Big]+
 \frac{N-1}{\rho^{\varepsilon}}\varepsilon^{-1}\rho^{1-\varepsilon}v'(\rho)\Big|^r d\rho
 \]
\[
 =\omega_{N-1}\int_0^1\varepsilon^{1-2r}\rho^{\varepsilon N-\frac{\varepsilon\beta}{q}-1}\Big|\rho^{2-2\varepsilon}
 \Big[v''(\rho)-(\varepsilon-1)\rho^{-1}v'(\rho)\Big]+
 (N-1)\varepsilon\rho^{1-2\varepsilon}v'(\rho)\Big|^r d\rho
 \]
 (here there was the mistake)
 \[
 =\omega_{N-1}\int_0^1\varepsilon^{1-2r}\rho^{\varepsilon N-\frac{\varepsilon\beta}{q}-1+2r-2r\varepsilon}
 \Big|v''(\rho)+\frac{N\varepsilon-2\varepsilon+1}{\rho}v'(\rho)\Big|^r d\rho
 \]
 \[
 =\omega_{N-1}\varepsilon^{1-2r}\int_0^1\rho^{\varepsilon N-\frac{\varepsilon\beta}{q}-1+2r-2r\varepsilon}
 \rho^{-(N\varepsilon-2\varepsilon+1)r}\Big|\Big({\rho}^{N\varepsilon-2\varepsilon+1}v'(\rho)\Big)'\Big|^r d\rho
 \]
 \[
 =\omega_{N-1}\varepsilon^{1-2r}\int_0^1\rho^{\varepsilon N-\frac{\varepsilon\beta}{q}+r-Nr-1}
 \left|\left({\rho}^{N\varepsilon-2\varepsilon+1}v'(\rho)\right)'\right|^r d\rho
 \]
 \[
 =\omega_{N-1}\varepsilon^{1-2r}\int_0^1\rho^{\varepsilon N-\frac{\varepsilon\beta}{q}+r-Nr\varepsilon-1}
 \Big|\Big({\rho}^{\gamma}v'(\rho)\Big)'\Big|^r d\rho.
 \]
 where $\gamma=N\varepsilon-2\varepsilon+1$. Moreover, by the choice of $\varepsilon$,
\[
 \int_{\Omega}|x|^{\alpha}|u(x)|^{p+1}\ dx=\omega_{N-1}\ \varepsilon\int_0^1|v(\rho)|^{p+1}\rho ^{N-1}d\rho.
 \]
Thus, we get the following estimate for the radial level:
 \begin{equation}
 S_{\a,\beta}^{rad}= C_N\varepsilon^{-2r-\frac{r}{p+1}+1}\inf_{v\in
 E_{rad}\setminus\{0\}}\frac{\int_0^1\rho^{\varepsilon N-\frac{\varepsilon\beta}{q}+r-Nr\varepsilon-1}
 \left|\left({\rho}^{\gamma}v'(\rho)\right)'\right|^r d\rho}{(\int_0^1|v(\rho)|^{p+1}\rho^{N-1})^{\frac{r}{p+1}}}
\end{equation}
It is now sufficient to show  that there exists $\eta>0$
independet on $\varepsilon$ such that
$$
\inf_{v\in
 E_{rad}\setminus\{0\}}\frac{\int_0^1\rho^{\varepsilon N-\frac{\varepsilon\beta}{q}+r-Nr-1}
 \left|\left({\rho}^{\gamma}v'(\rho)\right)'\right|^r
 d\rho}{(\int_0^1|v(\rho)|^{p+1}\rho^{N-1})^{\frac{r}{p+1}}}\ge\eta\varepsilon^{1/q}
 \qquad{\rm uniformly}\quad{\rm as}\quad \varepsilon\to 0
$$
We proceed as in the embedding result setting $w(\rho)=v'(\rho)\rho^{\gamma}$. Then
 \[
 |v(t)|=\Big|\int_{1}^tv'(\rho)d\rho\Big|=\Big|\int_{1}^tw(\rho)\rho^{-\gamma}d\rho\Big|=
 \Big|\int_{1}^t\rho^{-\gamma}\Big(\int_0^{\rho}w'(s)ds\Big)d\rho\Big|
 \]
 \[
 = \Big|\int_{1}^t\rho^{-\gamma}\Big[\int_0^{\rho}w'(s)
  s^{\frac{\varepsilon N-\frac{\varepsilon\beta}{q}+r-Nr\varepsilon-1}{r}}s^{-\frac{\varepsilon N-\frac{\varepsilon\beta}{q}+r-Nr\varepsilon-1}{r}}ds\Big]d\rho\Big|
 \]
 (H\"{o}lder inequality)
 \[
 \le \Big|\int_{1}^t\rho^{-\gamma}\Big(\int_0^{1}|w'(s)|^r s^{\varepsilon
 N-\frac{\varepsilon\beta}{q}+r-Nr\varepsilon-1}ds\Big)^{1/r}\Big(\int_0^{\rho}
s^{(-\varepsilon
N+\frac{\varepsilon\beta}{q}-r+Nr\varepsilon+1)q}ds\Big)^{\frac{1}{q+1}}d\rho\Big|
 \]
 \[
\le \Big|\int_{1}^t\rho^{-\gamma}\Big(\int_0^{1}|w'(s)|^r
s^{\varepsilon
 N-\frac{\varepsilon\beta}{q}+r-Nr\varepsilon-1}ds\Big)^{1/r}\Big(\int_0^{\rho}
s^{-\varepsilon Nq
+{\varepsilon\beta}-1+N(q+1)\varepsilon}ds\Big)^{\frac{1}{q+1}}d\rho\Big|
 \]
 \[
\le
\frac{1}{[\varepsilon(N+\beta)]^{\frac{1}{q+1}}}\Big|\int_{1}^t\rho^{-\gamma}\Big(\int_0^{1}|w'(s)|^r
s^{\varepsilon
N-\frac{\varepsilon\beta}{q}+r-Nr\varepsilon-1}ds\Big)^{1/r}
\rho^{-\frac{\varepsilon Nq}{q+1} +\frac{\varepsilon
\beta}{q+1}+N\varepsilon}d\rho\Big|
 \]
  \[
=\frac{1}{[\varepsilon(N+\beta)]^{\frac{1}{q+1}}}\Big(\int_0^{1}|w'(s)|^r
s^{\varepsilon
N-\frac{\varepsilon\beta}{q}+r-Nr\varepsilon-1}ds\Big)^{1/r}
 \Big|\int_{1}^t\rho^{2\varepsilon-\frac{Nq\varepsilon}{q+1}+\frac{\varepsilon\beta}{q+1}-1}d\rho\Big|
 \]
   \[
=\frac{1}{[\varepsilon(N+\beta)]^{\frac{1}{q+1}}}\Big(\int_0^{1}|w'(s)|^r
s^{\varepsilon
 N-\frac{\varepsilon\beta}{q}+r-Nr\varepsilon-1}ds\Big)^{1/r}
 \Big|\int_{1}^t\rho^{\varepsilon(2-N+\frac{N+\beta}{q+1})-1}d\rho\Big| =: \Im \vspace{0.3cm}
 \]
 For  ${q+1}\neq\frac{N+\beta}{ N-2}$ one has
 \[
\Im=\frac{1}{[\varepsilon(N+\beta)]^{\frac{1}{q+1}}}\Big(\int_0^{1}|w'(s)|^r
s^{\varepsilon
 N-\frac{\varepsilon\beta}{q}+r-Nr\varepsilon-1}ds\Big)^{1/r} \left|\frac{t^{{\varepsilon(2-N+\frac{N+\beta}{q+1})}}-1}
 {\varepsilon\Big(N-2-\frac{N+\beta}{q+1}\Big)}\right|
 \]
 Therefore
\[
\Big(\int_0^1|v(t)|^{p+1}t^{N-1}dt\Big)^{\frac{r}{p+1}}
\]
\[
\le\frac{1}{[\varepsilon(N+\beta)]^{\frac{1}{q}}}
\Big(\int_0^{1}|w'(s)|^r s^{\varepsilon
N-\frac{\varepsilon\beta}{q}+r-Nr\varepsilon-1}ds\Big)
 \bigg(\int_{0}^{1}\left|\frac{t^{{\varepsilon(2-N+\frac{N+\beta}{q+1})}}-1}
 {\varepsilon\Big(N-2-\frac{N+\beta}{q+1}\Big)}\right|^{p+1}
 t^{N-1}dt \bigg)^{\frac{r}{p+1}}
\]
Now we prove that the last term is uniformly bounded as
$\varepsilon\to 0 $. Let
\[
g_{\varepsilon}(t)=\Bigg|\frac{t^{{\varepsilon(2-N+\frac{N+\beta}{q+1})}}-1}
 {\varepsilon\left(N-2-\frac{N+\beta}{q+1}\right)}\Bigg|^{p+1}
 t^{N-1}
\]
We have that $$g_{\varepsilon}(t)\to(-\log t)^{p+1}t^{N-1}\ {\rm on}\ (0,1)\ , \ {\rm as}\ \varepsilon\to 0 \;$$
 and
 $$
 g_{\varepsilon}(t)\le(-\log t)^{p+1}\ {\rm on} \ (0,1) \ .
 $$
 Then, by the Dominated Convergence Theorem

$$
\int_0^1g_{\varepsilon}(t)dt\to\int_0^1(-\log t)^{p+1}t^{N-1}dt
$$
which is finite.

\end{proof}
\par \medskip \noindent



Quite surprisingly one can prove a similar result (with some
restrictions) even in the mixed case:

\begin{Theorem} (H\'enon-Hardy type system, M. Calanchi, B. Ruf (2010),\cite{CR}) \label{t3}
\par \noindent
Let $\alpha > 0$, $0 > \beta  > -N$, and suppose that $p , q > 1$.
\par \smallskip \noindent
a)  If $\frac{N+\alpha}{p+1}+\frac{N-|\beta|}{q+1} > N-2$ (i.e.
below the $\alpha\beta$-hyperbola),
\par\noindent then system (\ref{pb1}) has a radial solution (not necessarily of minimal energy)
\par \medskip \noindent
b) If $(p,q)$ satisfies \ $\frac{N}{p+1}+\frac{N -|\beta|}{q+1} >
N-2$, then $\inf_{E} R(u)$ is attained and hence system
(\ref{pb1}) has a ground state solution; furthermore, if  $p>\frac{N}{(N-1)q-1}$, and $\alpha >
0$ is sufficiently large, then the ground state solution is {\rm
not radially symmetric}.
\end{Theorem}

\begin{figure}[h]
 \centerline{\includegraphics[width=13cm]{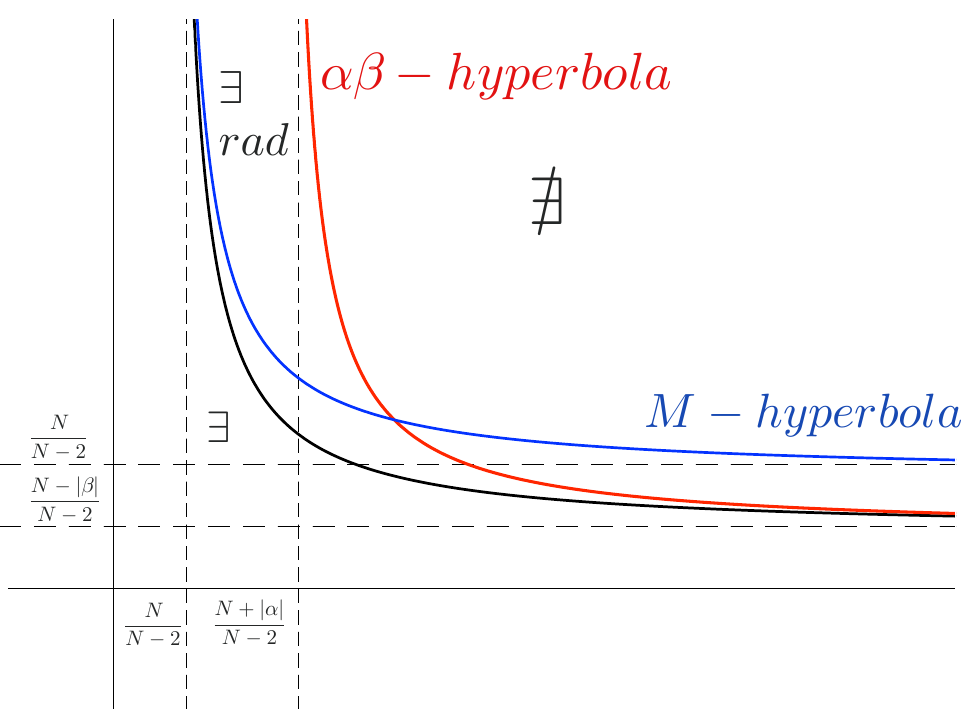}}
  \centerline{{\it Figure 3: Mixed case.}}
  \end{figure}

 The proof is similar of the one of the H\'enon
case. One only needs a compact embedding as in the Hardy case.

\par \medskip\noindent

\subsection{Open problems}

The present paper raises several open questions. A few list
follows:

{\bf{1.}} For systems of Hardy-H\'enon type:  suppose that
$\Omega$ is an arbitrary bounded domain: if $(p+1,q+1)$ lies below
the $\alpha\beta$-hyperbola and above the M-hyperbola (fig. 3),
then it is {\it not known} whether $\inf_E R(u)$ is attained.

2. For systems of H\'enon type:  in the subcritical case (under
the M-hyperbola) is it possible to describe the asymptotic
behavior of the ground state solution as $\alpha\to +\infty$?

3. For systems of H\'enon type:  is it possible to prove
multiplicity  of non-radial solutions  with some partial symmetry
in the spirit of Badiale and Serra's result (Theorem
\ref{badiale})? Which is in this case the new critical hyperbola?

 \vspace{.1cm}
and moreover

4. For the H\'enon equation in the supercritical case
$p<{\frac{N+1}{N-3}}$ (Theorem \ref{badiale}): which is the
asymptotic behavior of the partially symmetric solutions as
$\alpha\to +\infty$ or when $p\to{\frac{N+1}{N-3}}$?

 \vspace{2cm}

\newpage

%


\begin{thebibliography}{99}


\bibitem{AR} A. Ambrosetti, P. H. Rabinowitz,
 \textit{Dual Variational Methods in criticla point theory and applications}, J.  Funct. Anal.,
{\bf 14}, (1973), 349 - 381.


\bibitem{badialeserra} M.~Badiale, E.~Serra, \textit{Multiplicity
    results for the supercritical {H}\'enon equation}, Adv. Nonlinear
  Stud., \textbf{4}, (2004), No. 4, 453--467.
  
  
  
\bibitem{BDR} D. Bonheure, E. Dos Santos, M. Ramos, \textit{
  Symmetry and symmetry breaking for ground state solutions of some strongly coupled elliptic systems}, Journal of Functional Analysis \textbf{264}, Issue 1, (2013) 62--96.

  \bibitem{BN} H.~Brezis, L.~Nirenberg, \textit{Positive solutions of
  nonlinear elliptic equations involving critical Sobolev exponents},
  Comm. Pure Appl. Math, \textbf{36}, (1983),  437--477.


\bibitem{byeonwangI} J.~Byeon,  Z.-Q.~Wang,
\textit{On the {H}\'enon equation: asymptotic profile of ground
              states. {I}}, Ann. Inst. H. Poincar\'e Anal. Non Lin\'eaire, \textbf{23},
(2006), 803--828.

\bibitem{byeonwang} J.~Byeon, Z.-Q.~Wang, \textit{On the {H}\'enon
    equation: asymptotic profile of ground states. {II}}, J.
  Differential Equations, \textbf{216}, (2005), 78--108.


\bibitem{CR} M. ~Calanchi, B.~ Ruf,  \textit{Radial and non radial solutions
for Hardy-H\'enon type elliptic systems}, Calc. Var. Partial
Differential Equations, {\bf 38} (2010), no. 1-2, 111--133,

  \bibitem{CKN} L.~Caffarelli, R.~Kohn, L. Nirenberg \textit{First
order interpolation inequalities with weights}, Compositio
Math., \textbf{53}, (1984), 259--275.

\bibitem{caopeng} D.~Cao, S.~Peng, \textit{The asymptotic behaviour of
    the ground state solutions for {H}\'enon equation}, J. Math. Anal.
  Appl., \textbf{278}, (2003), 1--17.

  \bibitem{caopengyan} D.~Cao, S.~Peng, S.~Yan,
\textit{Asymptotic behaviour of ground state solutions for the H\'enon
equation}, IMA J. Appl. Math. \textbf{74} (2009), no. 3,
468--480.

\bibitem{Clem}
P.~Cl\'ement,  B.~de Pagter, G.~ Sweers, F.~ de Th\'elin, \textit{Existence
of solutions to a semilinear elliptic system through
Orlicz-Sobolev spaces}, Mediterr. J. Math.,  \textbf{1
}(2004), no. 3, 241--267.

\bibitem{deF-F} D.G. ~de Figueiredo, P. ~Felmer, \textit{On Superquadratic Elliptic Systems},
Trans. American Math. Soc., \textbf{343}, no. 1, (1994), 99--116.

\bibitem{deF-doO-R} D.G. ~de Figueiredo, J.M. ~do \'O, B. ~Ruf,
\textit{ An Orlicz-space approach to superlinear elliptic systems}, J.
Funct. Anal.  \textbf{224}  (2005),  no. 2, 471--496.

\bibitem{DPR} D.G. ~de Figueiredo, I. ~Peral, J. D. Rossi, \textit{The critical hyperbola
for a Hamiltonian system with weights},
 Ann. Mat. Pura Appl. (4)  \textbf{187}  (2008),  no.
3, 531--545.

\bibitem{eke}  I.~Ekeland, \textit{On the variational principle}, J. Math. Anal. Appl., \textbf{47}
(1974), 324--353.

\bibitem{GNN} B. Gidas, W.M. Ni and L. Nirenberg, \textit{Symmetry and Related Properties via the Maximum Principle},
Commun. Math. Phys. \textbf{ 68} (1979), 209--243.

\bibitem{HY} H.~ He, J. Yang, \textit{The Asymptotic Behavior of
Solutionsfor H\'enon systems with nearly critical exponent},
J. Math. Anal. Appl. \textbf{347} (2008), 459--471.

\bibitem{henon} M.~H\'{e}non, \textit{Numerical experiments on the
    stability of spheriocal stellar systems}, Astronomy and
  Astrophysics,  \textbf{24} (1973), 229--238.

\bibitem{H-vdV} J. ~Hulshof, R. Van der Vorst,
  \textit{Differential systems with strongly indefinite
  variational structure}, J. Funct. Anal.\textbf{114} (1993), 32--58.

 \bibitem{hulsvander2} J. ~Hulshof, R. Van der Vorst,
  \textit{Asymptotic behaviour of Ground States}, Proceed. Amer. Math. Soc. , \textbf{124} (1996), 2423--2431.

\bibitem{L} E.H. ~Lieb, \textit{Sharp Constants in the Hardy-Littlewood-Sobolev and Related Inequalities},
Annals of Mathematics, \textbf{118}, 1983, 349--374.

\bibitem{Lin} Lin, Chang Shou, \textit{ Interpolation inequalities with weights},
 Comm. Partial Differential Equations \textbf{11} (1986), no. 14, 1515--1538.

\bibitem{PLL} P-L.~Lions, \textit{The concentration-compactness
principle in the calculus of variations. The limit case. Part I},
Rev. Mat. Iberoamericana,  \textbf{1} (1985), 145--201.

\bibitem{liuyang}F. ~Liu, J. Yang, \textit{Nontrivial solutions of Hardy-Henon type elliptic systems}.
 Acta Math. Sci. Ser. B Engl. Ed., \textbf{27}  (2007),  no. 4, 673--688.

\bibitem{mitidieri} E. ~Mitidieri, \textit{A Rellich type identity and
applications},  Comm. Partial Differential Equations,
 \textbf{18}, (1993), 125--151.


\bibitem{ni} W.~M.~Ni, \textit{A nonlinear Dirichlet problem on the
unit ball and its applications},  Indiana Univ. Math. J.,
\textbf{31} (1982), 801--807.

\bibitem{opic} B.~Opic, A.~Kufner, \textit{Hardy-type inequalities}.  Pitman
      Research Notes in Mathematics Series, \textbf{219}, Longman
      Scientific \& Technical, 1990.
\bibitem{pacella}  F. Pacella, \textit{Symmetry results for solutions of semilinear elliptic equations with convex nonlinearities},
 J. Funct. Anal. \textbf{192} (2002), no. 1, 271--282.

\bibitem{pacella1} F. Pacella,  {\it Uniqueness of positive solutions of semilinear elliptic equations and related eigenvalue problems},
Milan J. Math. \textbf{73} (2005), 221--236.

\bibitem{PT}R. S.~Palais, C.L.~ Terng,  \textit{Critical point theory and submanifold geometry}, Lecture Notes in Mathematics, 1353. Springer-Verlag, Berlin, (1988).

\bibitem{peng} S.~Peng, \textit{Multiple boundary concentrating
    solutions to {D}irichlet problem of {H}\'enon equation}, Acta
  Math. Appl. Sin. Engl. Ser., \textbf{22}, (2006), 137--162.

\bibitem{serra03} E.~Serra, \textit{Non radial positive solutions for
    the {H}\'enon equation with critical growth}, Calc. Var. Partial
  Differential Equations, \textbf{23} (2005), No. 3, 301--326.

\bibitem{ssw} D.~Smets, J. ~Su, M.~Willem, \textit{Non-radial ground
    states for the H\'{e}non equation},  Communications in Contemporary
  Mathematics, \textbf{4}, no. 3 (2002), 467--480.

\bibitem{sw} D.~Smets, M.~Willem, \textit{Partial symmetry and
    asymptotic behavior for some elliptic variational problems}, Cal.
  Var. ,  \textbf{18} (2005), 57--75.

\bibitem{strauss} W.A.~Strauss, \textit{Existence of solitary waves in higher
      dimensions},  Comm. Math. Phys., \textbf{55}, (1977), 149--162.

\bibitem{st} M.~Struwe, \textit{Variational methods. Applications to nonlinear partial differential equations and Hamiltonian systems}. Third edition. Springer-Verlag, Berlin, (2000).

\bibitem{Xuan} B.~Xuan,  \textit{The solvability of
quasilinear Brezis-Nirenberg-type problems with singular weights}, Nonlinear Anal. TMA \textbf{62} (2005),\textbf{ 4}, 703--725.

\bibitem{wang} Xu-Jia ~Wang, \textit{Sharp constant in a Sobolev inequality},
Nonlinear Anal. TMA,    \textbf{20},  No. 3, (1993), 261--268.


\bibitem{ZH} Yajing  ~Zhang, Jianghao ~Hao \textit{The asymptotic behavior of the ground state solutions for  biharmonic equations,} Nonlinear Anal. TMA,   Vol.  \textbf{74} (2011),
2739--2749.

\end{thebibliography}
\end{document}